\documentclass[11pt]{amsart}

\usepackage{geometry}
\usepackage{amsmath,amssymb,bm}
\usepackage{booktabs}   
\usepackage{graphicx}

\theoremstyle{definition}
\newtheorem{theorem}{Theorem}[section]

\newtheorem{lemma}[theorem]{Lemma}

\newtheorem{remark}[theorem]{Remark}

\usepackage{color}

\newcommand{\T}{\mathcal{T}}
\newcommand{\E}{\mathcal{E}}
\newcommand{\intel}[2][\T_h]{\left(#2\right)_{#1}}
\newcommand{\intfa}[2][\partial\T_h]{\langle #2\rangle_{#1}}

\newcommand{\jump}[1]{[\![#1]\!]}
\newcommand{\dd}{\partial}

\newcommand{\qq}{\bm{q}}
\newcommand{\qqhat}{\widehat{\qq}_h}
\newcommand{\uhat}{\widehat{u}_h}

\newcommand{\vv}{\bm{v}}
\newcommand{\nn}{\bm{n}}
\newcommand{\rr}{\bm{r}}

\newcommand{\PN}{\bm{P}_{\bm N}}
\newcommand{\PV}{\bm{P}_{\bm V}}
\newcommand{\PW}{P_{W}}
\newcommand{\Eq}{\bm e_{\bm q}}
\newcommand{\Eqh}{\bm e_{\widehat{\bm q}}}
\newcommand{\Eu}{e_u}

\begin{document}
\title[A flux-based HDG method]{A flux-based HDG method}

\author{Issei Oikawa}
\address{Department of Mathematics, Institute of Pure and Applied Sciences, University of Tsukuba,
  1-1-1 Tennodai, Tsukuba, Ibaraki 305-8571, Japan}
\email{ioikawa00@gmail.com}


\maketitle
\begin{abstract}
  In this paper, we present a flux-based formulation of the  hybridizable discontinuous Galerkin (HDG) method
   for steady-state diffusion problems and propose a new method derived by letting a stabilization parameter tend to infinity.  
  Assuming an inf-sup condition, we prove its well-posedness and error estimates of optimal order.
  We show that the inf-sup condition is satisfied by some triangular elements. 
  Numerical results are also provided to support our theoretical results.
\end{abstract}

\section{Introduction}
We consider the hybridizable discontinuous Galerkin (HDG) method for the steady-state diffusion problem with Dirichlet
boundary condition 
\begin{subequations}\label{prob}
  \begin{align}
    \qq + \nabla u   & = 0   \qquad \text{ in } \Omega,       \\
    \nabla \cdot \qq & = f    \qquad\text{ in } \Omega,       \\
    u               & = 0  \qquad\text{ on } \partial\Omega,
  \end{align}
\end{subequations}
where $\Omega \subset \mathbb{R}^d~(d=2,3)$ is a bounded convex polygonal or polyhedral
domain and $f$  is a given function.
In the original HDG method \cite{CGL2009},
a numerical trace $\widehat u_h$ is
introduced as an unknown variable to approximate the trace of $u$ on element boundaries, which corresponds to 
a Dirichlet boundary condition, and 
a numerical flux $\qqhat$ is properly defined.
The other variables $u_h$ and $\qq_h$ approximating  $u$ and $\qq$, respectively, can be eliminated
in element-by-element fashion and
we obtain a globally-coupled system of equations only in terms of $\widehat{u}_h$, which is called
static condensation.

In \cite{Cockburn2016}, a flux-based formulation is presented, in which 
the trace of $\qq$ on element boundaries instead of $\widehat{u}_h$ 
is hybridized, in other words, $\qqhat$ is an unknown variable and $\widehat{u}_h$ is defined 
in terms of $\qqhat$ and the other variables.
The flux-based method is a rewrite of 
the original HDG method and provides the same solution,
however, the local problem has a Neumann boundary condition, so that the static condensation is different from that of the original method. We note that the local solvability of the flux-based method is not obvious, which 
will be verified in Section \ref{subsect-localsol}.

In this paper, we propose a new flux-based method derived 
by passing the stabilization parameter to infinity.
In our method, $\qqhat$ is unknown and the numerical trace is defined by $\widehat{u}_h = u_h$. Since our  method has saddle point structure, 
its well-posedness depends on whether the inf-sup condition we define in Section \ref{sect-fluxhdg} is satisfied.
The inf-sup condition is fulfilled 
 if we use triangular meshes and the polynomials of degree $k$ and $k+1$ for $\qqhat$ and $u_h$, respectively, with a non-negative integer $k$. 
In addition, the proposed method using such approximation spaces achieves the optimal convergence rates in $u_h$ and $\qq_h$, like the HDG method with the so-called Lehrenfeld--Sch\"{o}berl stabilization 
proposed in \cite{Lehrenfeld2010} and analyzed in \cite{Oikawa2015,Oikawa2018}.
Although, in the Lehrenfeld--Sch\"{o}berl stabilization, the $L^2$-orthogonal projection onto the approximation space of $\widehat{u}_h$ is inserted in front of $u_h$, 
such a projection is not used in our method because it is naturally incorporated through the transmission condition in a flux-based formulation. 

The rest of the paper is organized as follows.
In Section 2, we introduce notation and present the flux-based formulation, and a new method is derived from it. We verify the local solvability of the methods.
In Section 3, we establish a priori estimate and error estimates of optimal order for our method, assuming that an inf-sup condition holds. 
In Section 4, we prove that the inf-sup condition is satisfied if the polynomial degrees for $\qqhat$ and $u_h$ are $k$ and $k+1$, respectively, and triangular meshes are used.
In Section 5, numerical results are presented to validate our theoretical results.
%
%
\section{A flux-based HDG formulation}
\subsection{Notation}
To begin with, we introduce notation to present the HDG method via flux hybridization.
Let $\{\T_h\}_h$ be a family of meshes satisfying the quasi-uniform condition, where $h$
stands for the mesh size.
Let $\E_h$ denote the set of all edges or faces of elements in $\T_h$.
Let $L^2(\E_h)$ denote the $L^2$-space on $\bigcup_{e \in \E_h} e$ and $P_m(\T_h)$ and $P_m(\E_h)$ denote the spaces of element-wise and edge-wise polynomials of degree $m$, respectively.

We use the usual symbols of Sobolev spaces \cite{AdFo2003}, such as
$H^m(D)$, $H^m(D)^d$, $\|\cdot\|_{m, D} := \|\cdot\|_{H^m(D)}$, and
$|\cdot|_{m, D} := |\cdot|_{H^m(D)}$ for a domain $D$ and an integer $m$.
We may omit the subscripts when $D = \Omega$ or $m=0$, such as $\|\cdot\|_{m} = \|\cdot\|_{m,\Omega}$,
$\|\cdot\| = \|\cdot\|_{0,\Omega}$, and $|\cdot|_{m} = |\cdot|_{m,\Omega}$.
The piecewise Sobolev space of order $m$ is denoted by $H^m(\T_h)$.
The inner products are defined as
\begin{align*}
& (\qq,\vv)_K  = \int_K \qq\cdot\vv dx, \quad
   \intel{\bm q, \bm v} = \sum_{K\in \T_h} (\bm q,\vv)_K, \\
& (u,w)_K  = \int_K uwdx, \quad
   \intel{u,w} = \sum_{K\in\T_h} (u,w)_K,  \\
& \intfa[\dd K]{u,w} = \int_{\dd K} uw ds, \quad
   \intfa{u,w} = \sum_{K \in \T_h} \intfa[\dd K]{u,w}, \\
& \intfa[\E_h]{\mu,\lambda} = \sum_{e \in \E_h} \int_e \mu\lambda ds, \quad
    (u,w) = \int_\Omega uwdx.
\end{align*}
We define the induced norms from these inner products by
\[
  \|\vv\|_{\T_h} = \intel{\vv,\vv}^{1/2},\quad
     \|w\|_{\T_h} = \intel{w,w}^{1/2},\quad 
  \|w\|_{\dd\T_h} = \intfa{w,w}^{1/2}.
\]
Throughout the paper, we use the symbol $C$ to denote
a generic constant independent of the mesh size $h$ 
and $\nn$ to stand for the unit (outer) normal vector to an edge $e \in \E_h$ or 
 $\dd K$ for $K \in \T_h$.
%
%
\subsection{Finite element spaces}
Let $k$ be a non-negative integer.
We define the local approximation spaces on $K \in \T_h$ as 
\begin{align*}
\bm{V}(K) = P_k(K)^d, \quad 
W(K) = P_{k+1}(K),
\end{align*}
where $P_m(K)$ stands for the space of polynomials of degree $m$.
We introduce an approximation space  for  $\qq|_e \cdot \nn$ on 
$e \in \E_h$, 
\[ 
\bm{N}(e) = P_k(e),
\]
and assume that $(I -\nn\otimes \nn) \rr = \bm 0$ for $\rr \in \bm{N}(e)$.
The global finite element spaces are defined by
\begin{align*}
& \bm V_h :=  \{ \vv \in  L^2(\Omega)^d :  \vv|_K \in \bm  V(K) ~ \forall K \in \T_h \}, \\
& W_h := \{w \in L^2(\Omega) :  w|_K \in W(K) ~ \forall K \in \T_h\}, \\
& \bm N_h := \{\rr \in L^2(\E_h)^d :   \rr  |_{e} \in \bm  N(e) ~ \forall e \in \E_h\}.
\end{align*}

Let $\PV$, $\PW$, and $\PN$ denote the $L^2$-projections onto $\bm V_h$, $W_h$, and $\bm{N}_h$, respectively. 
The following approximation properties hold
for $1\le s \le k+1$: If $u \in H^{k+2}(\Omega)$, then
\begin{subequations} \label{apppr}
  \begin{align}
    \|\qq - \PV \qq\| & \le Ch^s |\qq|_{s},
    \label{apppr-a} \\
    \|\qq\cdot \nn - (\PV\qq)\cdot\nn\|_{\dd\T_h} & \le Ch^{s-1/2}|\qq|_{s},
    \label{apppr-b} \\
    \|u - \PW u \| & \le Ch^{s}|u|_{s}, 
    \label{apppr-c} \\
    \|u - \PW u \|_{\dd\T_h} & \le Ch^{s+1/2}|u|_{s+1}, 
    \label{apppr-d} \\
    \|(\qq- \PN\qq)\cdot\nn \|_{\dd\T_h} & 
    \le Ch^{s+1/2}|\qq|_{s+1}. 
    \label{apppr-e}
  \end{align}
\end{subequations}
\subsection{A flux-based HDG method} \label{sect-fluxhdg}
The solution of the original HDG method, 
$(\qq_h, u_h, \widehat u_h) \in \bm V_h \times W_h  \times M_h$,
is defined by
\begin{subequations} \label{hdg-u}
\begin{align}
\intel{\qq_h,\vv} - \intel{u_h, \nabla \cdot \vv}
+ \intfa{\widehat u_h, \vv \cdot \nn}
& = 0               
&& \forall \vv \in \bm V_h,
    \label{hdg-u-v} \\
    -\intel{\qq_h, \nabla w} +\intfa{\qqhat\cdot \nn, w}
& = (f,w)          
&& \forall w \in W_h,     \label{hdg-u-w}\\
    \intfa{\qqhat\cdot\nn, \mu}& =0
&& \forall \mu \in M_h,  \label{hdg-u-mu} \\
\qqhat\cdot\nn := \qq_h\cdot\nn + \tau(u_h - \widehat{u}_h) 
& \text{ on } \dd K
&& \forall K \in \T_h,  \label{hdg-qflux}
\end{align}
\end{subequations}
where $M_h$ is an approximation space for the trace $u|_{\E_h}$ 
and $\tau$ is a stabilization parameter.
Another formulation via flux hybridization is also stated in \cite{Cockburn2016}, 
which reads as follows: Find
$(\qq_h, u_h, \widehat \qq_h) \in \bm V_h \times W_h \times \bm N_h$ such that
\begin{subequations} \label{hdg-q}
  \begin{align}
    \intel{\qq_h,\vv} - \intel{u_h, \nabla \cdot \vv}
    + \intfa{\uhat, \vv \cdot \nn}
& = 0 && \forall \vv \in \bm V_h,
    \label{hdg-q-v} \\
    -\intel{\qq_h, \nabla w} +\intfa{\widehat{\qq}_h\cdot \nn, w}
& = (f,w)    && \forall w \in W_h,
    \label{hdg-q-w} \\
    \intfa{\uhat, \rr\cdot\nn} & =0
&& \forall \rr \in \bm N_h, \label{hdg-q-mu} \\
    \uhat :=  u_h + \tau^{-1} (\qq_h - \widehat{\qq}_h)\cdot\nn
& \text{ on } \dd K 
&& \forall K \in \T_h.\label{hdg-q-uflux}
\end{align}
\end{subequations}
This method is a rewrite of the original HDG method and its solution coincides with that of the original method.
We can verify that by expressing the hybrid variables in terms of 
$u_h$ and $\qq_h$. Let $K^+$ and $K^-$ be adjacent elements 
sharing an internal edge $e \in \E_h$ and let $\nn^+$ and $\nn^-$ denote 
the outer unit normal vectors to $\dd K^+$ and $\dd K^-$, respectively.
For a function $w$, let $w^+$ and $w^-$ stand for the trace of $(w|_{K^+})|_e$ and $(w|_{K^-})|_e$, respectively.
In both methods, $\uhat$ and $\qqhat$ are single valued on element boundaries
from the transmission conditions \eqref{hdg-u-mu} and \eqref{hdg-q-mu}.
From \eqref{hdg-qflux} or  \eqref{hdg-q-uflux}, it follows that 
\begin{align*}
  \qqhat\cdot\nn^+ &= \qq_h^+\cdot\nn^+ + \tau(u_h^+ - \uhat),\\
  \qqhat\cdot\nn^- &= \qq_h^-\cdot\nn^- + \tau(u_h^- - \uhat).
\end{align*}
Solving these equations, we have 
\begin{align*}
  \uhat &= \frac{1}{2}(u_h^+ + u_h^-) + \frac{1}{2\tau}(\qq_h^+\cdot\nn^{+} + \qq_h^-\cdot\nn^-), \\ 
  \qqhat\cdot\nn^\pm &=  \frac{1}{2}(\qq_h^+ + \qq_h^-)\cdot\nn^\pm
  + \frac{\tau}{2}(u_h^+ \nn^+ + u_h^-\nn^-)\cdot\nn^\pm.
\end{align*}
Therefore, we see that the equations \eqref{hdg-u-v}-\eqref{hdg-u-w} and \eqref{hdg-q-v}-\eqref{hdg-q-w}
give the same solution $u_h$ and $\qq_h$. 
However, note that the procedures of the static condensation are different and the local solvability of the flux-based method is not obvious, which we will prove later.

We consider the limiting case of $\tau \to +\infty$ in \eqref{hdg-q}. 
In this case, \eqref{hdg-q-uflux} is naturally interpreted as $\uhat = u_h$, 
which leads to the following scheme:  
Find $(\qq_h, u_h, \qqhat) \in \bm{V}_h\times W_h \times \bm{N}_h$ such that
\begin{subequations} \label{hdg-qinf}
  \begin{align}
    \intel{\qq_h,\vv} - \intel{u_h, \nabla \cdot \vv}
    + \intfa{u_h, \vv \cdot \nn}
  & = 0 && \forall \vv \in \bm V_h,
    \label{hdg-qinf-v}\\
    -\intel{\qq_h, \nabla w} +\intfa{\widehat{\qq}_h\cdot \nn, w}
& = (f,w)           && \forall w \in W_h,
    \label{hdg-qinf-w} \\
    \intfa{u_h, \rr\cdot\nn} & =0
&& \forall \rr \in \bm N_h.  \label{hdg-qinf-mu}                                 
\end{align}
\end{subequations}
We remark that the above method is not always well-posed. 
Assume that $f \equiv 0$. 
By taking $\vv = \qq_h$,
$w = u_h$, and $\rr = \qqhat$ in \eqref{hdg-qinf}, 
we have  $\qq_h = \bm{0}$. 
From \eqref{hdg-qinf-v} with $\qq_h = \bm{0}$ and 
\eqref{hdg-qinf-mu}, it follows that $u_h = 0$.
However, $\qqhat$  still remains 
unknown, which depends on if the following equation implies $\qqhat\cdot\nn = 0$:
\begin{align} \label{qnw}
\intfa{\widehat{\qq}_h\cdot \nn, w} &= 0 
 \qquad \forall w \in W_h.
\end{align}
Indeed, when $\bm N_h$ and $W_h$ are piecewise constant spaces, it is easy to see that there exists $\qqhat \in \bm{N}_h$ satisfying \eqref{qnw} and  
$\qqhat\cdot\nn \neq 0$.
For this reason, we need the following inf-sup condition for the well-posedness: There 
exists a constant $C$ independent of $h$ such that, for all $\rr \in \bm{N}_h$, 
\begin{align} \label{infsup-rw}
\|h^{1/2}\rr\cdot\nn\|_{\dd\T_h} 
  \le C \sum_{K\in\T_h} \sup_{w\in W(K)} \frac{\intfa[\dd K]{\rr\cdot\nn, w}}{\|\nabla w\|_{L^2(K)} + \|h^{-1/2}\PN w\|_{L^2(\dd K)}},   
 \end{align}
where $\PN w := \PN(w|_{\dd K} \nn)\cdot\nn$. 
In order to derive a priori estimates, we also use the transposed version of the inf-sup condition: There exists a constant $C$ such that  
\begin{align} \label{infsup-wr}
\|h^{-1/2}\PN w\|_{\dd\T_h} 
  \le C \sum_{K\in\T_h} \sup_{\rr\in \bm{N}(\dd K)} \frac{\intfa[\dd K]{\rr\cdot\nn, w}}{\|h^{1/2}\rr\cdot\nn\|_{L^2(\dd K)}} \qquad 
  \forall w \in W_h.
 \end{align}

\subsection{Local solvability} \label{subsect-localsol}
We here verify the local solvability 
of the flux-based methods \eqref{hdg-q} and \eqref{hdg-qinf}.

We first consider the local problem of \eqref{hdg-q}.
Let $\widehat{\qq}_{\dd K}$ denote the restriction of $\widehat{\qq}_h$ to $\dd K$
and let $\jump{w}$ denote the jump of a function $w$. We define
$\|\mu\|_{\E_h} = \intfa[\E_h]{\mu,\mu}^{1/2}$
and
$\|\rr\|_{\E_h} = \intfa[\E_h]{\rr,\rr}^{1/2}$.
We introduce the mean-zero subspace of $W(K)$,
\[
  W_0(K) := \{ w \in W(K) : (w,1)_K = 0 \}.
\]
The local problem reads:
Find $(\qq_K, u_{K0}) \in \bm V(K)\times W_0(K)$ such that
\begin{subequations} \label{localprob-naumann}
  \begin{align}
    (\qq_K,\vv)_K
    + \intfa[\dd K]{\tau^{-1}\qq_{K}\cdot\nn,\vv\cdot\nn} + (\nabla u_{K0}, \vv)_K
     & = \intfa[\dd K]{\tau^{-1}\widehat{\qq}_{\dd K}\cdot\nn,\vv\cdot\nn}
     & \forall \vv \in \bm V(K),
    \label{localprob-naumann-v}
    \\
    -(\qq_K, \nabla w_0)_K
     & = (f,w_0)_K
    - \intfa[\dd K]{\widehat{\qq}_{\dd K}\cdot\nn,w_0}
     & \forall w_0 \in  W_0(K).
    \label{localprob-naumann-w}
  \end{align}
\end{subequations}
To verify the well-posedness of the local problem, 
we let $f \equiv 0$ and $\widehat{\bm{q}}_{\dd K}\cdot\nn = 0$.
Taking  $\vv = \qq_K$ in \eqref{localprob-naumann-v} and $w_0 = u_{K0}$ in \eqref{localprob-naumann-w}, we have $\qq_K = \bm{0}$.
Since we can $\vv = \nabla u_{K0}$, we get $\nabla u_{K0} = 0$, 
which implies $u_{K0} = 0$.
Therefore, $\qq_K$ and $u_{K0}$ are uniquely determined if $\widehat \qq_{\dd K}$ is given.

However, the piecewise constant part of $u_h$, denoted by $\overline{u}_h$, remains unknown.
Eliminating $\qq_K$ and $u_{K0}$ in element-by-element fashion by static condensation,
we obtain the global equations for 
$(\widehat{\qq}_h, \overline{u}_h) \in \bm N_h \times P_0(\T_h)$
\begin{subequations} \label{globalprob}
  \begin{align}
    \intfa{\overline{u}_h - \tau^{-1}\widehat{\qq}_h \cdot\nn, \rr\cdot\nn}
     & =-\intfa{u_{h0} - \tau^{-1}\qq_h \cdot\nn, \rr\cdot\nn}
    =: F_1(\rr)
     &                                                         & \forall \rr \in \bm N_h,
    \label{globalprog-neumann-r}                                                                     \\
    \intfa{\widehat{\qq}_h\cdot\nn,\overline{w}}
     & =
    (f,\overline{w})
     &                                                         & \forall \overline{w} \in P_0(\T_h),
    \label{globalprog-neumann-wbar}
  \end{align}
\end{subequations}
where $u_{h0} := u_h - \overline{u}_h$ and 
we note that $u_{h0}$ has been determined by \eqref{localprob-naumann}.
We will show that the global problem is well posed.
To this end, we first prove the following inf-sup condition.
\begin{theorem} \label{infsup-w0jump}
  There exists a
  positive constant $C$ such that
  \[
    C \|\jump{\overline{w}}\|_{\E_h} \le
    \sup_{\rr\in\bm N_h}\frac{\intfa{\rr\cdot\nn, \overline{w}}}{\|\rr\cdot\nn\|_{\dd\T_h}} 
    \qquad \forall \overline{w} \in P_0(\E_h).
  \]
\end{theorem}
\begin{proof}
  For $\overline{w} \in P_0(\T_h)$ with $\overline{w} \neq 0$, we define $\rr = \jump{\overline{w}}$.
  Then, it follows that
  \begin{align*}
    \intfa{\rr\cdot\nn, \overline{w}} =
    \intfa[\E_h]{\rr, \jump{\overline{w}}} =
    \|\jump{\overline{w}}\|_{\E_h}^2 =
    \|\rr\|_{\E_h}\|\jump{\overline{w}}\|_{\E_h}.
  \end{align*}
  Since $\|\rr\|_{\E_h} \le \|\rr\cdot\nn\|_{\dd\T_h} \le \sqrt{2}\|\rr\|_{\E_h}$
   for $\rr \in \bm N_h$, we have
  \begin{align*}
    \|\jump{\overline{w}}\|_{\E_h}
    \le  \frac{\intfa{\rr\cdot\nn, \overline{w}}}{\|\rr\|_{\E_h}}
    \le \sqrt{2}\cdot \frac{\intfa{\rr\cdot\nn, \overline{w}}}{\|\rr\cdot\nn\|_{\dd\T_h}},
  \end{align*}
  which completes the proof.
\end{proof}
We now prove an a priori estimate for the global problem \eqref{globalprob} using the inf-sup condition 
in Theorem \ref{infsup-w0jump},
which ensures that  the problem admits a unique solution.
\begin{theorem} \label{esti-neumann}
  There exists a positive constant $C$ independent of $h$ such that
  \begin{align*}
    \|\tau^{-1/2}\widehat{\qq}_h\cdot\nn\|_{\dd \T_h}
    +\|\jump{\overline{u}_h}\|_{\E_h}
    \le C\left( \|{u}_{h0}\|_{\dd\T_h} +  \|\qq_h\cdot\nn\|_{\dd\T_h}+\|f\|\right).
  \end{align*}
\end{theorem}
\begin{proof}
Taking $\rr = \widehat{\qq}_h$ in \eqref{globalprog-neumann-r} and
$\overline{w} = \overline{u}_h$ in \eqref{globalprog-neumann-wbar}, we have
\begin{align*}
\intfa{\overline{u}_h, \widehat{\qq}_h\cdot\nn}
- \|\tau^{-1/2}\widehat{\qq}_h\cdot\nn\|_{\dd \T_h}^2
 & = F_1(\widehat{\qq}_h),
\\
\intfa{\widehat{\qq}_h\cdot\nn, \overline{u}_h}
 & = (f,\overline{u}_h).
\end{align*}
It then follows that
\begin{align*}
\|\tau^{-1/2}\widehat{\qq}_h\cdot\nn\|_{\dd \T_h}^2
 & =  (f, \overline{u}_h)    - F_1(\widehat{\qq}_h)
\le \|f\|\|\overline{u}_h\|
+ \|F_1\| \|\widehat{\qq}_h\cdot\nn\|_{\dd \T_h},
\end{align*}
where
  \[
    \|F_1\| := \sup_{\rr\in\bm N_h}\frac{F(\rr)}{\|\rr\cdot\nn\|_{\dd\T_h}}.
  \]
Using Young's inequality, we deduce
\begin{align} \label{ineq-qn-u}
\|\tau^{-1/2}\widehat{\qq}_h\cdot\nn\|_{\dd \T_h}^2
\le C\left( \|F_1\|^2 + \epsilon^{-2}\|f\|^2 \right) 
  + \epsilon^2\|\overline{u}_h\|_{\dd\T_h}^2
\end{align}
for any $\epsilon>0$.
By Theorem \ref{infsup-w0jump}, \eqref{globalprog-neumann-r} and \eqref{ineq-qn-u}, we have
\begin{align*}
\|\jump{\overline{u}_h}\|_{\E_h}
 & \le C \sup_{\rr\in\bm N_h}\frac{\intfa{\rr\cdot\nn, \overline{u}_h}}{\|\rr\cdot\nn\|_{\dd\T_h}} \\
 & \le C \sup_{\rr\in\bm N_h}
\frac{F_1(\rr) + \intfa{\tau^{-1}\widehat{\qq}_h\cdot\nn,\rr\cdot\nn}}{\|\rr\cdot\nn\|_{\dd\T_h}}  \\
 & \le C(\|F_1\| + \|\tau^{-1}\widehat{\qq}_h\cdot\nn\|_{\dd\T_h})                                 \\
 & \le C\left(\|F_1\| + \epsilon^{-1} \|f\| + \epsilon\|\overline{u}_h\|_{\dd\T_h}\right).
\end{align*}
Since both $\|\jump{\,\cdot\,}\|_{\E_h}$ 
and $\|\cdot\|_{\dd\T_h}$ are norms on $P_0(\T_h)$, 
they are equivalent to each other and $\|\overline{u}_h\|_{\dd\T_h}$ is bounded by $C\|\jump{\overline{u}_h}\|_{\E_h}$. Choosing $\epsilon$ sufficiently small, 
we get 
    \[
      \|\jump{\overline{u}_h}\|_{\E_h}
      \le C\left(\|F_1\| + \|f\|\right).
    \]
From this and \eqref{ineq-qn-u}, it follows that
  \[
    \|\tau^{-1/2}\widehat{\qq}_h\cdot\nn\|_{\dd \T_h}
    \le C\left(\|F_1\| +  \|f\|\right).
  \]
  Since we can bound as $\|F_1\| \le C(\|u_{h0}\|_{\dd\T_h} + \|\qq_h\cdot\nn\|_{\dd\T_h})$, the proof is complete.
\end{proof}

Next, we show that the local problem of the proposed method is well-posed.
Let us define $\bm{R}_0(\dd K) = \{ w|_{\dd K} \nn: w \in P_0(K)\}$, whose dimension is one.
The local problem of \eqref{hdg-qinf} is as follows:
Find $(\qq_K, u_{K0}, \overline{u}_K) \in \bm V(K)\times W_0(K)\times P_0(K)$ such that
 \begin{subequations} \label{localprob-qinf}
  \begin{align}
    (\qq_K,\vv)_K
     + (\nabla u_{K0}, \vv)_K
     & = 0
     && \forall \vv \in \bm V(K),
     \label{localprob-qinf-v}
    \\
    -(\qq_K, \nabla w_0)_K
     & = (f,w_0)_K
    - \intfa[\dd K]{\widehat{\qq}_{\dd K}\cdot\nn,w_0}
     && \forall w_0 \in  W_0(K), 
     \label{localprob-qinf-w} \\
     \intfa[\dd K]{\overline{u}_K +u_{K0}, \rr_0 \cdot\nn}
    &= 0
    && \forall \rr_0 \in \bm{R}_0(\dd K).
     \label{localprob-qinf-r}
  \end{align}
\end{subequations}
The well-posedness is verified by setting all terms on the right-hand side to zero and a straightforward computation.
\section{Error Analysis}
\subsection{A priori estimate} We consider the proposed method in general form 
\begin{subequations} \label{hdg-qinf-gen}
  \begin{align}
    \intel{\qq_h,\vv} - \intel{u_h, \nabla \cdot \vv}
    + \intfa{u_h, \vv \cdot \nn}
                             & = F_1(\vv)   &                          & \forall \vv \in \bm V_h,
    \label{hdg-qinf-v-gen}                                                                                    \\
    -\intel{\qq_h, \nabla w} +\intfa{\widehat{\qq}_h\cdot \nn, w}
                             & = F_2(w)           &                          & \forall w \in W_h,
    \label{hdg-qinf-w-gen}                                                                                    \\
    \intfa{u_h, \rr\cdot\nn} & = F_3(\rr)
                             &                   & \forall \rr \in \bm N_h,
    \label{hdg-qinf-mu-gen}                                                        
    \end{align}
\end{subequations}
where $F_1 : \bm{V}_h \to \mathbb{R}$, 
$F_2 : W_h \to \mathbb{R}$, and $F_3 : \bm{N}_h \to \mathbb{R}$ 
are linear functionals and their norms are defined by
\begin{align*}
  \|F_1\| &= \sup_{\vv\in \bm V_h} \frac{F_1(\vv)}{\|\vv\|}, \\ 
  \|F_2\| &= \sum_{K \in \T_h}\sup_{w\in W(K)} 
  \frac{F_2(w)}{\|\nabla w\|_{L^2(K)} 
  + \|h^{-1/2}\PN w\|_{L^2(\dd K)}}, \\ 
  \|F_3\| &= \sup_{\rr\in \bm N_h} \frac{F_3(\rr)}{\|h^{1/2}\rr\cdot\nn\|_{\dd\T_h}}.
\end{align*}
We first establish a priori estimate for the problem.
\begin{theorem} \label{thm-stab}
  Let $(\qq_h, u_h, \qqhat) \in \bm V_h \times W_h \times \bm N_h$ be 
  a solution of \eqref{hdg-qinf-gen}.  
  Then there exists a constant $C$ such that
  \begin{align*}
    \|\qq_h\| +\|h^{1/2}\qqhat\cdot\nn\|_{\dd\T_h} 
    + \|\nabla u_h\|_{\T_h} + \|h^{-1/2}\PN u_h\|_{\dd\T_h}
     &\le C\left(\|F_1\|+\|F_2\| + \|F_3\|\right).
  \end{align*}
\end{theorem}
\begin{proof}
From the inf-sup condition \eqref{infsup-rw} and \eqref{hdg-qinf-w-gen}, it follows that 
\begin{equation} \label{esti-qhatn}
\begin{aligned}
  \|h^{1/2}\qqhat\cdot\nn\|_{\dd\T_h} 
  &\le C \sum_{K \in \T_h} 
  \sup_{w \in W(K)} \frac{\intfa[\dd K]{\qqhat\cdot\nn, w}}{\|\nabla w\|_{L^2(K)} + \|h^{-1/2}\PN w\|_{L^2(\dd K)}} \\ 
  &= C \sum_{K \in \T_h} 
  \sup_{w \in W(K)} \frac{F_2(w) + \intel[K]{\qq_h, \nabla w}}{\|\nabla w\|_{L^2(K)} + \|h^{-1/2}\PN w\|_{L^2(\dd K)}}\\
  &\le C\left( \|F_2\| + \|\qq_h\|\right).
\end{aligned}
\end{equation}
Integrating by parts in \eqref{hdg-qinf-v-gen}, we have
\begin{align}
\intel{\qq_h,\vv}
 + \intel{\nabla u_h, \vv} =  F_1(\vv).
 \label{hdg-qinf-v-gen2}
\end{align}
Substituting $\vv = \nabla u_h$ in the above equation, we get 
\begin{equation} \label{grad-uh}
\|\nabla u_h\|_{\T_h} \le C  \left( \|\qq_h\| + \|F_1\|\right).
\end{equation}
Taking $\vv = \qq_h$ in \eqref{hdg-qinf-v-gen2}, 
$\rr = \qqhat$ in \eqref{hdg-qinf-mu-gen}, and $w = u_h$ in \eqref{hdg-qinf-w-gen}, we have 
\begin{align*}
 \|\qq_h\|^2 
 &= F_1(\qq_h) + F_2(u_h) - F_3(\qqhat) \\
 &\le \|F_1\|\|\qq_h\| 
    + C\|F_2\| \sum_{K\in\T_h} 
    \left(\|\nabla u_h\|_{L^2(K)} 
    + \|h^{-1/2} \PN u_h\|_{\dd K}\right) 
    + \|F_3\|\|h^{1/2}\qqhat\cdot\nn\|_{\dd\T_h}.
\end{align*}
By \eqref{infsup-wr}, we have
\begin{align*}
\|h^{-1/2} \PN u_h\|_{\dd\T_h}
 & \le C \sum_{K\in\T_h} \sup_{\rr \in \bm{N}_h}
 \frac{\intfa{\rr\cdot\nn, \PN u_h}}{\|h^{1/2}\rr\cdot\nn\|} \\ 
 &=  C \sum_{K\in\T_h} \sup_{\rr \in \bm{N}_h}
 \frac{\intfa{\rr\cdot\nn, u_h}}{\|h^{1/2}\rr\cdot\nn\|} \\ 
 &\le C \sum_{K\in\T_h} \sup_{\rr \in \bm{N}_h}
 \frac{F_3(\rr)}{\|h^{1/2}\rr\cdot\nn\|} \\ 
 &\le C\|F_3\|.
\end{align*}
Combining this with \eqref{grad-uh}, we get
\[
 \sum_{K\in\T_h} 
    \left(\|\nabla u_h\|_{L^2(K)} 
    + \|h^{-1/2} \PN u_h\|_{\dd K}\right) 
    \le C \left(\|\qq_h\| + \|F_3\|\right).
\]
Thus we estimate as 
\begin{align*}
 \|\qq_h\|^2 & \le \|F_1\|\|\qq_h\| + C\|F_2\| \left(\|\qq_h\| + \|F_1\|\right)
  + C\|F_3\|(\|F_2\|+\|\qq_h\|). 
\end{align*}
Using Young's inequality, we obtain 
\begin{align*}
\|\qq_h\|^2 \le C \left( \|F_1\|^2 + \|F_2\|^2 + \|F_3\|^2\right). 
\end{align*}
From this and  \eqref{esti-qhatn}, it follows that 
\begin{align*}
  \|h^{1/2}\qqhat\cdot\nn\|_{\dd\T_h} \le C\left( 
    \|F_2\|+  \|\qq_h\|\right) \le C(\|F_1\|+\|F_2\|+\|F_3\|),
\end{align*}
which completes the proof.
\end{proof}
If $f\equiv 0$, it follows from Theorem \ref{thm-stab}
that $\qq_h = \bm{0}$,  $u_h$ is constant on each element, 
 $\widehat{\qq}_h\cdot\nn = 0$, and $\PN u_h = 0$ on element boundaries. 
Thus, we have verified the existence and
uniqueness of our method.
\subsection{Optimal convergence of $\qq_h$}
The projections of errors are defined as
\begin{align*}
  \Eq = \PV\qq-\qq_h, \quad
  \Eu = \PW u - u_h, \quad
  \Eqh\cdot\nn = \PN\qq\cdot\nn - \widehat{\qq}_h\cdot\nn.
\end{align*}
\begin{theorem} \label{thm-errest-dq}
If $u \in H^{k+2}(\Omega)$, then  we have
  \begin{align*}
    \|\Eq\| +   \|h^{1/2}\Eqh\cdot\nn\|_{\dd\T_h}&\le Ch^{k+1}|u|_{k+2}.
  \end{align*}
\end{theorem}
\begin{proof}
  The problem \eqref{prob} is rewritten into
  \begin{align*}
    \intel{\qq,\vv} - \intel{u, \nabla\cdot \vv} + \intfa{u ,\vv\cdot\nn} 
    & = 0
    &         & \forall \vv \in \bm V_h, \\
    -\intel{\qq, \nabla w} + \intfa{\qq\cdot\nn, w}                       
    & = (f,w)
    &         & \forall w \in W_h,       \\
    \intfa{u, \rr\cdot\nn} & = 0
    &         & \forall \rr \in \bm N_h.
  \end{align*}
By the property of the $L^2$-projections, the above equations become
\begin{subequations}\label{eqs-hdgproj}
\begin{align}
\intel{\PV\qq,\vv} 
-\intel{\PW u, \nabla \cdot \vv} + \intfa{\PW u ,\vv\cdot\nn} 
  & =G_1(\vv)
   && \forall \vv \in \bm V_h, \\
-\intel{\PV\qq, \nabla w} + \intfa{\PN \qq\cdot\nn, w}                               
& = (f,w) + G_2(w)
   && \forall w \in W_h,       \\
\intfa{\PW u, \rr\cdot\nn}                                                        
  & = G_3(\rr)
   && \forall \rr \in \bm N_h,     
\end{align}
\end{subequations}
where we have integrated by parts in the first equation and 
\begin{align*}
    G_1(\vv) &:= - \intfa{u - \PW u, \vv\cdot\nn}, \\ 
    G_2(w) &:= - \intfa{\qq\cdot\nn - \PN\qq\cdot\nn, w}, \\ 
    G_3(\rr) &:= - \intfa{u - \PW u, \rr\cdot\nn}. 
\end{align*}
The norms of $G_1$ and $G_3$ are bounded as 
\begin{align*} 
    \|G_1\| &= \sup_{\vv\in\bm V_h} \frac{G_1(\vv)}{\|\vv\|} 
    \le \|u - \PW u\|_{\dd\T_h} \cdot Ch^{-1/2}
    \le Ch^{k+1}|u|_{k+2}, 
    \\
    \|G_3\| &= \sup_{\rr\in\bm N_h} \frac{G_3(\rr)}{\|h^{1/2}\rr\cdot\nn\|_{\dd\T_h}} 
    \le \|u - \PW u\|_{\dd\T_h} \cdot h^{-1/2}
    \le Ch^{k+1}|u|_{k+2}. 
\end{align*}
Using \cite[Lemma 3]{Oikawa2015}, we can estimate 
\begin{align*}
|G_2(w)| &\le |\intfa{(\qq - \PN\qq)\cdot\nn, w}| \\  
        &=\|(\qq - \PN\qq)\cdot\nn\|_{\dd\T_h} \cdot Ch^{1/2}\|\nabla w\|_{\T_h} .
\end{align*}
The norm of $G_2$ is bounded as
\begin{align*}
\|G_2\| &= \sum_{K\in\T_h} \sup_{w\in W(K)} \frac{G_2(w)}{\|\nabla w\|_{L^2(K)}+\|h^{-1/2}\PN w\|_{L^2(\dd K)}} 
\le Ch^{k+1}|\qq|_{k+1}. 
\end{align*}
Subtracting \eqref{hdg-qinf} from \eqref{eqs-hdgproj}, we obtain the error equations
\begin{subequations}
\label{erroreqs}
\begin{align}
\intel{\Eq,\vv} 
  -\intel{\Eu, \nabla\cdot\vv} + \intfa{\Eu,\vv\cdot\nn} 
  & = G_1(\vv)
   && \forall \vv \in \bm V_h,
  \label{erroreqs-v}  \\
  -\intel{\Eq, \nabla w} + \intfa{\Eqh\cdot\nn, w}                  
  & =  G_2(w)
  && \forall w \in W_h, \label{erroreqs-w} \\
\intfa{\Eu,\rr\cdot\nn}  &= G_3(\rr)          
    && \forall\rr \in \bm N_h.
       \label{erroreqs-r} 
\end{align}
\end{subequations}
Applying Theorem \ref{thm-stab} to the error equations leads to 
\begin{align*}
    \|\Eq\| + \|h^{1/2}\Eqh\cdot\nn\|_{\dd\T_h} 
    &\le C\left(\|G_1\|+\|G_2\|+\|G_3\| \right) 
     \le Ch^{k+1}|u|_{k+2}.
\end{align*}
\end{proof}
From Theorem 3.1, it also follows that 
\begin{align} \label{erresti-gradu}
\|\nabla \Eu\|_{\T_h} + \|h^{-1/2}\PN\Eu\|_{\dd\T_h}\le C h^{k+1}|u|_{k+2}.
\end{align}
\subsection{$L^2$-error estimate of $u_h$}
We consider the following adjoint problem:
Find $(\bm\theta, \xi) \in \bm H^1(\Omega)\times (H^2(\Omega)\cap H_0^1(\Omega))$ such that
\begin{subequations}
  \begin{align*}
    \bm\theta + \nabla \xi & = \bm 0  ~\quad \text{ in } \Omega, \\
    \nabla \cdot \bm\theta & = \Eu   \quad \text{ in } \Omega,   \\
    \xi                    & = 0 ~\quad \text{ on } \dd\Omega.
  \end{align*}
\end{subequations}
It is well known that the elliptic regularity holds:
\[
  \|\bm\theta\|_{1} + \|\xi\|_2 \le C\|\Eu\|.
\]
We provide an $L^2$-error estimate of $u_h$ by the Aubin--Nitsche technique.
\begin{theorem} \label{err-u}
  If $k \ge 1$ and $u \in H^{k+2}(\Omega)$, then there exists a constant $C$ such that 
  \begin{align*}
    \|\Eu \| \le Ch^{k+2}|u|_{k+2}.
  \end{align*}
\end{theorem}
\begin{proof}
Since \eqref{eqs-hdgproj} holds for the adjoint problem, we have
\begin{subequations}
\label{eqs-hdgproj-adj}
\begin{align}
  \intel{\PV\bm\theta,\vv} 
  - \intel{\PW \xi, \nabla\cdot \vv}
   + \intfa{\PW \xi ,\vv\cdot\nn} & 
   = G'_1(\vv)
    && \forall \vv \in \bm V_h, \label{eqs-adj-v}\\
  -\intel{\PV\bm\theta, \nabla w} 
  + \intfa{\PN\bm\theta\cdot\nn, w}                        
  & = (\Eu,w) + G'_2(w)
    && \forall w \in W_h, \label{eqs-adj-w}\\
      \intfa{\PW \xi, \rr\cdot\nn}          
  & =   G_3'(\rr)
   && \forall \rr\in\bm N_h,
\label{eqs-adj-r}
\end{align}
\end{subequations}
where 
\begin{align*}
G_1'(\vv) &:= - \intfa{\xi - \PW\xi, \vv\cdot\nn}, \\ 
G_2'(w) &:= - \intfa{(\bm\theta - \PN\bm\theta)\cdot\nn, w}, \\ 
G_3'(\rr) &:= - \intfa{\xi - \PW \xi, \rr\cdot\nn}.  
\end{align*}
Taking $\vv = -\Eq$ in \eqref{eqs-adj-v},
 $w = \Eu$ in \eqref{eqs-adj-w}, and $\rr = \Eqh$ in \eqref{eqs-adj-r},
we have
\begin{equation}
\label{adj-1}
\begin{aligned}
  -\intel{\PV\bm\theta, \Eq} 
  - \intel{\nabla\PW\xi, \Eq}
  -\intel{\PV\bm\theta,\nabla\Eu}
  + \intfa{\PN\bm\theta\cdot\nn, \Eu} 
  +\intfa{\PW\xi, \Eqh\cdot\nn}\\ 
  =\|\Eu\|^2 - G_1'(\Eq) + G_2'(\Eu) + G_3'(\Eqh).
\end{aligned}
\end{equation}
Choosing $\vv = -\PV\bm\theta$ in \eqref{erroreqs-v},
$w = \PW\xi$ in \eqref{erroreqs-w}, and $\rr = \PN\bm\theta$ in \eqref{erroreqs-r},
we have
\begin{equation} \label{adj-2}
\begin{aligned}
  -\intel{\Eq, \PV\bm\theta}
  -\intel{\nabla\Eu, \PV\bm\theta}
  - \intel{ \Eq, \nabla\PW\xi} 
  + \intfa{\Eqh\cdot\nn,  \PW\xi} 
  + \intfa{\Eu,\PN\bm\theta\cdot\nn}\\
= - G_1(\PV\bm\theta) + G_2(\PW\xi) + G_3(\PN\bm\theta).
\end{aligned}
\end{equation}
Subtracting \eqref{adj-2} from \eqref{adj-1} yields
\begin{align*}
\|\Eu\|^2 & =
 - G_1(\PV\bm\theta) + G_2(\PW\xi) + G_3(\PN\bm\theta) 
 - G_1'(\Eq) + G_2'(\Eu) + G_3'(\Eqh).
\end{align*}
We will bound the terms on the right-hand side. 
The first and third terms are bounded as
\begin{align*}
|- G_1(\PV\bm\theta) + G_3(\PN\bm\theta)|
&= |\intfa{u - \PW u, (\PV \bm\theta - \PN\bm\theta)\cdot\nn}| \\ 
&\le  \|u - \PW u\|_{\dd\T_h} (\|(\PV \bm\theta -\bm\theta)\cdot\nn\|_{\dd\T_h} 
    + \|(\bm\theta- \PN\bm\theta)\cdot\nn\|_{\dd\T_h}) \\
&\le Ch^{k+3/2}|u|_{k+2} \cdot Ch^{1/2}|\bm\theta|_1 \\
&= Ch^{k+2} |u|_{k+2}|\bm\theta|_1.
\end{align*}
Let $P_1$ denote the $L^2$-projection from $L^2(\Omega)$ onto $P_1(\T_h)$. 
Note that  
$
  \intfa{(\qq - \PN\qq)\cdot\nn, P_1\xi} = 0
$
since we assume $k \ge 1$.
We have
\begin{equation} \label{G2PWxi}
\begin{aligned}
|G_2(\PW\xi)| 
&= |\intfa{(\qq - \PN\qq)\cdot\nn, \PW\xi - P_1\xi}| \\ 
&\le \|(\qq - \PN\qq)\cdot\nn\|_{\dd\T_h}
   (\|\PW\xi -\xi\|_{\dd\T_h} + \|\xi - P_1\xi\|_{\dd\T_h}) \\ 
&\le Ch^{k+1/2}|\qq|_{k+1} \cdot Ch^{3/2}|\xi|_{2} \\ 
& = Ch^{k+2}|\qq|_{k+1}|\xi|_{2}.
\end{aligned}
\end{equation}
The rest terms are 
bounded as follows:
  \begin{align*}
  |G_1'(\Eq)| 
  &\le Ch |\xi|_2 \|\Eq\|, \\ 
  |G_2'(\Eu)|  
  &\le C h |\bm\theta|_1 \|\nabla \Eu\|_{\T_h}  
  \qquad \text{ (by \cite[Lemma 3]{Oikawa2015})}\\ 
  &\le C h |\bm\theta|_1 \left( \|\Eq\|+ h^{k+1}|u|_{k+2}\right),
    \qquad \text{(by \eqref{erresti-gradu})} \\ 
  |G_3'(\Eqh)| 
  &\le C h |\xi|_2 \|h^{1/2}\Eqh\cdot\nn\|_{\dd\T_h}. 
\end{align*}
Thus we deduce 
  \[
   \|\Eu\| \le C\left(h \|\Eq\| + h \|h^{1/2}\Eqh\cdot\nn\|_{\dd\T_h} + h^{k+2}|u|_{k+2}\right)  
   \le C h^{k+2}|u|_{k+2}.
  \]
\end{proof}
\begin{remark}
Theorem \ref{err-u} also holds for $k=0$. 
When $k=0$, we can use the Crouzeix--Raviart interpolation $I_{CR}$ instead of $\PW$. Then, 
the right-hand sides of \eqref{erroreqs} and \eqref{eqs-hdgproj-adj} are 
changed as 
$G_1(\vv) = -\intel{\nabla (u - I_{CR}u), \vv}$, 
$G_3(\rr) = 0$, $G'_1(\vv) = -\intel{\nabla(\xi - I_{CR}\xi), \vv}$, $G_3'(\rr) = 0$. 
It is clear that they are bounded by the Schwarz inequality and the interpolation error estimate, and 
we do not need to use the $L^2$-projection 
in \eqref{G2PWxi} 
since $G_2(I_{CR}\xi) = \intfa{(\qq-\PN\qq)\cdot\nn, I_{CR}\xi} = 0$.
\end{remark}
\section{Proof of the inf-sup condition for triangular elements}
We show that the inf-sup conditions \eqref{infsup-rw} and \eqref{infsup-wr} are satisfied for the triangular $P_k$-element in the two-dimensional case.

\begin{figure}[h]
    \centering
    \includegraphics{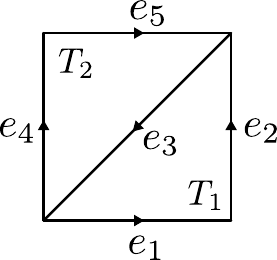}
    \caption{The reference triangles $T_1$ and $T_2$ and the orientation of the edges $e_1, \ldots, e_5$}
    \label{fig-elem}
\end{figure}

Let $T_1$ and $T_2$ be the reference triangles 
whose vertices are $\{(0,0), (0,1), (1,1)\}$ and $\{(1,0),(0,0),(0,1)\}$, respectively.
Let $\{e_1, e_2, e_3\}$ and $\{e_3, e_4, e_5\}$ denote be the edges of $T_1$
and $T_2$, respectively, see \figurename \ \ref{fig-elem}.
For $1 \le i \le 5$, we define $F_i$ by the 
 linear transforms from $e_i$ to $[-1,1]$ such that
 \begin{align*} 
 F_1((0,0)) = -1, \quad F_1((1,0)) = 1, \\ 
 F_2((1,0)) = -1, \quad F_2((1,1)) = 1, \\ 
 F_3((1,1)) = -1, \quad F_3((0,0)) = 1, \\ 
 F_4((0,0)) = -1, \quad F_4((0,1)) = 1, \\ 
 F_5((0,1)) = -1, \quad F_5((1,1)) = 1.
 \end{align*}

Let $\chi^{(i)}$ denote the characteristic function of $e_i$ and 
 let $\varphi_p$ denote the Legendre polynomial of degree $p$ on $[-1,1]$.
We define the normalized Legendre polynomial on $e_i$ by 
\[ 
\varphi_p^{(i)} = \frac{(\varphi_p \circ F_i) 
    \chi^{(i)}}{\left(\int_{e_i}(\varphi_p \circ F_i)^2 ds\right)^{1/2}}
\]
for $0 \le p \le k$.
 It is clear that $\{\varphi_0^{(i)}, \varphi_1^{(i)}, \ldots, \varphi_k^{(i)}\}$ is a basis of $\bm{N}(e_i) = P_k(e_i)$ and satisfies the orthogonality
\[
  \intfa[\dd T_1\cup\dd T_2]{\varphi_p^{(i)}, \varphi_q^{(j)}} = 
  \delta_{ij}\delta_{pq} \quad \text{ for } 1 \le i,j,p,q \le 5. 
\]
Note that $\varphi_p^{(i)} = (-1)^p$ at the starting point of $e_i$ and $\varphi_p^{(i)} = 1$ at the end point of $e_i$.

\begin{lemma} \label{lem-infsup-r} 
Let $\rr \in  \bigoplus_{1 \le i \le 5} \bm{N}(e_i)$. It holds that $\rr\cdot\nn = 0$ on $e_i$
 $(1 \le i \le 5)$ 
 if and only if 
 \begin{equation} \label{infsup-eq1}
     \intfa[\dd T_1 \cup \dd T_2]{\rr\cdot \nn, w} = 0 
 \quad \forall w \in W(T_1)\oplus W(T_2).
 \end{equation}
 
\end{lemma}

\begin{proof}
We prove that only when $k$ is even, i.e., $k = 2k'$, since the proof when $k$ is odd is similar.
We show that $\rr\cdot\nn = 0$ follows from \eqref{infsup-eq1}. We can write $\rr\cdot\nn$ as 
\begin{align*}
\rr\cdot \nn = \sum_{i=1}^5\sum_{p=0}^{2k'} a_p^{(i)}\varphi_p^{(i)}, \quad a_p^{(i)} \in \mathbb{R}.
\end{align*}
First, we show that $a_{2q}^{(i)} = 0$ for $1\le q \le k'$ and $1 \le i \le 5$.
Let us define  
\[ 
w_1 = \varphi_{2q}^{(1)}
-\varphi_{2k'+1}^{(2)} 
+ \varphi_{2k'+1}^{(3)}
\] 
for $1 \le q \le k'$, see also \figurename\ \ref{fig-testfunction}.
Since $w_1$ is continuous at the vertices, 
there exists $Ew_1 \in W(T_1) \oplus W(T_2)$ 
such that $Ew_1 = w_1$ on element boundaries. 
For simplicity, we use the same symbol $w_1$ to denote $Ew_1$.  
Choosing $w = w_1$ in \eqref{infsup-eq1} and noting that $\intfa[\dd T_1\cup\dd T_2]{\rr\cdot\nn, \varphi_{2k'+1}^{(i)}} = 0$ for any $i$, 
we have
\begin{align*}
    \intfa[\dd T_1 \cup \dd T_2]{\rr\cdot \nn, w_1} 
    &=
    \sum_{p=0}^{2k'} \intfa[e_1]{ a_p^{(1)}\varphi_p^{(1)}, \varphi_{2q}^{(1)}}
     = a_{2q}^{(1)} = 0.
\end{align*}
Taking $w = \varphi_{2k'+1}^{(1)} 
+\varphi_{2q}^{(2)}
-\varphi_{2k'+1}^{(3)}$ 
and 
$w = -\varphi_{2k'+1}^{(1)}
+ \varphi_{2k'+1}^{(2)}
+ \varphi_{2q}^{(3)}$,
we get $a_{2q}^{(2)} = a_{2q}^{(3)} = 0$ for $1 \le q \le k'$.
Similarly, it follows that 
$a_{2q}^{(4)} = a_{2q}^{(5)} = 0$ for $1 \le q \le k'$.

\begin{figure}[bht]
    \centering
    \includegraphics[width=400pt]{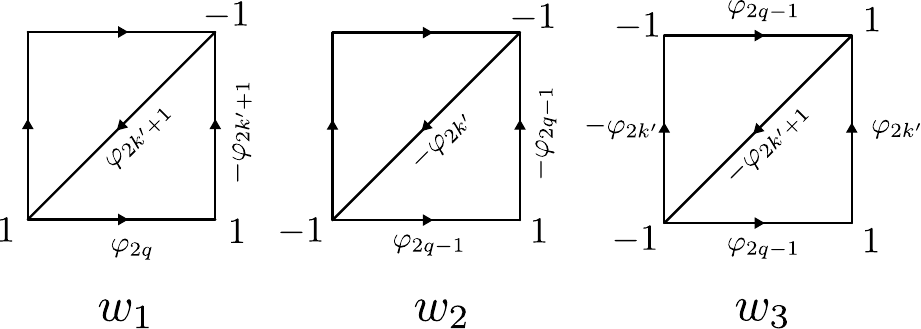}
    \caption{Diagrams of the test functions $w_1$, $w_2$, and $w_3$, where the subscripts are discarded and the number at each vertex indicates the value of the function at the vertex}
    \label{fig-testfunction}
\end{figure}

Next, we show that $a_{1}^{(i)} = a_{3}^{(i)} = \cdots = a_{2k'-1}^{(i)}$ for $1 \le i \le 5$.
In \eqref{infsup-eq1}, choosing 
\[ 
w = w_2 := \varphi_{2q-1}^{(1)}
-\varphi_{2q-1}^{(2)} + \varphi_{2k'}^{(3)},
\]
in view of $a_{2k'}^{(i)} = 0$ for $1 \le i \le 5$, we have
\begin{align*}
\intfa[\dd T_1 \cup \dd T_2]{\rr\cdot \nn, w_2}  
= a_{2q-1}^{(1)} - a_{2q-1}^{(2)} = 0    \quad (1 \le q \le k').
\end{align*}
Choosing $w = \varphi_{2k'}^{(1)} + \varphi_{2q-1}^{(2)} -\varphi_{2q-1}^{(3)}$ 
  in \eqref{infsup-eq1} yields $a_{2q-1}^{(2)} = a_{2q-1}^{(3)}$ for $1 \le q \le k'$. 
Similarly, we deduce that $a_{2q-1}^{(3)} = a_{2q-1}^{(4)} = a_{2q-1}^{(5)}$ for $1 \le q \le k'$.
Hence, we define $a_{2q-1} := a_{2q-1}^{(1)} = \cdots = a_{2q-1}^{(5)}$ for $1 \le q \le k'$
 to omit the subscripts.
 
Finally, we prove that $a_{2q-1} = 0$ for $1 \le q \le k'$.
We take the following $w_3$ as a test function: 
\[
 w_3 := \varphi_{2q-1}^{(1)} + \varphi_{2k'}^{(2)}
  -\varphi_{2k'+1}^{(3)} -\varphi_{2k'}^{(3)}
  + \varphi_{2q-1}^{(5)}. 
\]
Since $w_3$ is single valued on $e_3$, we see that 
\[ 
  \intfa[\dd T_1 \cup \dd T_2]{\rr\cdot\nn, w_3|_{e_3}} =
  \intfa[\dd T_1 \cup \dd T_2]{\rr\cdot\nn, -\varphi_{2k'+1}^{(3)}} = 0.
\]
Noting that $a_{2k'}^{(i)} = 0$ for $1 \le i \le 5$, we get
\[ 
\intfa[\dd T_1 \cup \dd T_2]{\rr\cdot\nn, \varphi_{2k'}^{(2)}}
=\intfa[\dd T_1 \cup \dd T_2]{\rr\cdot\nn, 
-\varphi_{2k'}^{(4)}} = 0
\]
and 
\begin{align*}
    \intfa[\dd T_1 \cup \dd T_2]{\rr\cdot \nn, w_3} 
    &=
    \intfa[\dd T_1 \cup \dd T_2]{\rr\cdot \nn, \varphi_{2q-1}^{(1)} + \varphi_{2q-1}^{(5)}} 
    =
    2 a_{2q-1} = 0.
\end{align*}
Thus, we conclude that all coefficients equal zero.
\end{proof}

\begin{lemma} \label{lem-infsup-local}
Let $K_1$ and $K_2 \in \T_h$ be two adjacent triangles. 
There exists a constant $C$ independent $h$ such that, 
for $\rr \in \bigoplus_{e \in \dd K_1 \cup \dd K_2} \bm{N}(e)$,
\[
  \sum_{i=1}^2 \|h^{1/2}\rr\cdot\nn\|_{L^2(\dd K_i)} 
  \le C\sum_{i=1,2} \sup_{w\in W(K_i)} \frac{\intfa[\dd K_i]{\rr\cdot\nn, w}}{\|\nabla w\|_{L^2(K_i)} + \|h^{-1/2}\PN \nabla w\|_{L^2(\dd K_i)}}. 
\]
\end{lemma}
\begin{proof}
Let $T_1$ and $T_2$ be the reference triangles and 
let $e_1, e_2,\ldots, e_5$ be the edges of $T_1$ and $T_2$. By Lemma \ref{lem-infsup-r}, we see that 
\[
   \|\rr\cdot\nn\|'_{\dd T_1 \cup \dd T_2} := 
   \sum_{i=1,2}
   \sup_{w \in W(T_i)}  
   \frac{\intfa[\dd T_i]{\rr\cdot\nn, w}}{\|\nabla w\|_{L^2(T_i)} + \|\PN w\|_{\dd T_i}}
\]
is a norm on $\bigoplus_{1\le i \le 5} \bm{N}(e_i)$.
Since any two norms on a finite-dimensional space are equivalent, 
there exists a constant $C$ such that
\[
  \|\rr\cdot\nn\|_{\dd T_1 \cup \dd T_2} 
  \le C \|\rr\cdot\nn\|'_{\dd T_1 \cup \dd T_2} 
  \quad \forall \rr\in \bigoplus_{1 \le i \le 5} \bm{N}(e_i).
\]
By considering the Piola transforms from $T_i$ to $K_i$ $(i=1,2)$ and the scaling argument, we obtain the assertion. 
\end{proof}
The inf-sup condition \eqref{infsup-rw} immediately follows from Lemma \ref{lem-infsup-local}. Similarly, we can prove the transposed inf-sup condition \eqref{infsup-wr} from the following lemma.

\begin{lemma}
Let $K$ be an element of $\T_h$ and $w \in W(K)$.
Then, $\PN w = 0$ on $\dd K$ if and only if  
\[
   \intfa[\dd K]{\rr\cdot\nn, w} = 0 \qquad \forall \rr \in \bm{N}(\dd K).
\]
\end{lemma}
\begin{proof}
Since $\PN w \in \bm{N}(\dd K)$, we can choose $\rr \cdot \nn = \PN w$. 
We then have 
\begin{align*}
    \intfa[\dd K]{\rr\cdot\nn, w} 
    = \intfa[\dd K]{\rr\cdot\nn, \PN w} 
    = \|\PN w\|_{\dd K}^2 = 0,
\end{align*}
which implies $\PN w|_{\dd K} = 0$.
\end{proof}
The proof of the transposed inf-sup condition is the same as in Lemma \ref{lem-infsup-r}, 
so we omit it here.
\begin{remark}
We have proved that the inf-sup condition \eqref{infsup-rw} holds for the triangular elements where $\bm{N}_h$ and $W_h$ are polynomials of degree $k$ and $k+1$, respectively. 
However, it is still open whether there exists a pair of $\bm{N}_h$ and $W_h$ satisfying the inf-sup condition 
\eqref{infsup-rw} in the three- or higher-dimensional cases.
\end{remark}
\section{Numerical results}
In this section, we examine the convergence property
of the proposed method \eqref{hdg-qinf} by numerical experiments.
As a test problem, we consider the Poisson equation with  homogeneous Dirichlet boundary condition 
\begin{align*}
 -\Delta u &= 2\pi^2\sin(\pi x)\sin(\pi y) && \text{ in } \Omega, \\
         u &= 0 && \text{ on } \partial\Omega,
\end{align*}
where $\Omega = (0,1)^2$ and the exact solution is given by $u(x,y) = \sin(\pi x)\sin(\pi y)$.
We use the unstructured triangulations whose mesh sizes are approximately $0.4\times 2^{-i}$ $(1 \le i \le 4)$ and 
compute the solution of \eqref{hdg-qinf}, varying the polynomial degree $k$  from $0$ to $3$.
All numerical computations are carried out by FreeFEM \cite{freefem}.
The $L^2$-errors of $\qq_h$ and $u_h$  are displayed in \tablename \ \ref{tb-hist}. 
We observe that the orders of convergence in $\qq_h$ and $u_h$ are $k+1$ and $k+2$, respectively, 
which are of optimal order and fully agrees with Theorems \ref{thm-errest-dq} and \ref{err-u}.

\begin{table}[hb]
\centering
\label{tb-hist}
\caption{Convergence history for the method \eqref{hdg-qinf}}
\begin{tabular}{cccccc}  \toprule
$k$ &	$h$	&	$\|\qq - \qq_h\|$	&	Order	&	$\|u - u_h\|$	& Order	\\ 
\midrule
	&	0.1901	&	3.472E-01	& ---	&	1.027E-02	&	---	\\
0 &	0.1025	&	1.660E-01	&	1.20	&	2.404E-03	&	2.35	\\
	&	0.0509	&	8.355E-02	&	0.98	&	6.150E-04	&	1.94	\\
	&	0.0262	&	4.120E-02	&	1.07	&	1.461E-04	&	2.17	\\ 
\midrule
   &	0.1901	&	1.827E-02	& ---	&	2.868E-04	&	---	\\
1 &	0.1025	&	5.218E-03	&	2.03	&	4.389E-05	&	3.04	\\
	&	0.0509	&	1.221E-03	&	2.07	&	4.926E-06	&	3.12	\\
	&	0.0262	&	2.967E-04	&	2.14	&	5.872E-07	&	3.21	\\
\midrule
    &	0.1901	&	1.172E-03	& ---	&	1.537E-05	&	---	\\
2 &	0.1025	&	1.212E-04	&	3.68	&	7.598E-07	&	4.87	\\
	&	0.0509	&	1.465E-05	&	3.01	&	4.524E-08	&	4.02	\\
	&	0.0262	&	1.797E-06	&	3.17	&	2.732E-09	&	4.24	\\
\midrule												
  &	0.1901	&	3.283E-05	&	---	&	2.773E-07	&	---	\\
3 &	0.1025	&	2.662E-06	&	4.07	&	1.185E-08	&	5.11	\\
  &	0.0509	&	1.408E-07	&	4.19	&	3.021E-10	&	5.23	\\
  &	0.0262	&	8.376E-09	&	4.27	&	1.326E-11	&	4.73	\\ 
  \bottomrule
\end{tabular}
\end{table}
\section*{acknowledgements}
The author would like to thank the anonymous reviewers for providing valuable comments and
suggestions.
%

\bibliographystyle{siam}
\bibliography{references}
\end{document}